\RequirePackage{etoolbox}
\csdef{input@path}{{style/}{graphics/}}
\documentclass[aap,MSNbibl,dvips]{arximspdf}
\makeatletter
   \@ifpackageloaded{graphicx}{}{\usepackage{graphicx}}
\makeatother

\doi{10.1214/14-AAP1024} 
\volume{25}
\issue{3}
\pubyear{2015}
\firstpage{1325}
\lastpage{1348}

\makeatletter
\renewcommand{\mid}{|}
\newcommand{\nicefrac}[2]{{#1}/{#2}}
\newtheorem{lm}{Lemma}

\newtheorem{prop}{Proposition}
\newproclaim{rem}{Remark}
\newproclaim{nota}{Notation}
\makeatother

\begin{document}
\begin{frontmatter}

\title{The internal branch lengths of\\ the Kingman
coalescent\thanksref{t1}}
\runtitle{Internal lengths}

\begin{aug}
\author[A]{\fnms{Iulia}~\snm{Dahmer}\corref{}\ead[label=e1]{dahmer@math.uni-frankfurt.de}\thanksref{T2}}
\and
\author[A]{\fnms{G\"otz}~\snm{Kersting}\ead[label=e2]{kersting@math.uni-frankfurt.de}}
\runauthor{I. Dahmer and G. Kersting}
\affiliation{Goethe-Universit\"at Frankfurt}
\address[A]{Institut f\"ur Mathematik\\
Goethe-Universit\"at Frankfurt\\
Robert-Mayer-Str. 10\\
Box 187\\
D-60325 Frankfurt am Main\\
Germany\\
\printead{e1}\\
\phantom{E-mail: }\printead*{e2}}
\end{aug}
\thankstext{T1}{Supported in part by the DFG Priority Programme SPP 1590 ``Probabilistic Structures in Evolution.''}
\thankstext{T2}{Supported by the German Academic Exchange Service (DAAD).}

\received{\smonth{4} \syear{2013}}
\revised{\smonth{12} \syear{2013}}

%
\begin{abstract}
In the Kingman coalescent tree the length of order $r$ is defined as
the sum of the lengths of all branches that support $r$ leaves. For
$r=1$ these branches are external, while for $r \ge2$ they are
internal and carry a subtree with $r$ leaves. In this paper we prove
that for any $s \in\mathbb N$ the vector of rescaled lengths of orders
$1\le r \le s$ converges to the multivariate standard normal
distribution as the number of leaves of the Kingman coalescent tends to
infinity. To this end we use a coupling argument which shows that for
any $r\ge2$ the (internal) length of order $r$ behaves asymptotically
in the same way as the length of order 1 (i.e., the external length).
\end{abstract}

%

\end{frontmatter}

\section{Introduction and main result}\label{sec1}

The Kingman coalescent was introduced in~\cite{Ki82} as a model for
describing the genealogical relationships between the individuals for a
wide class of population models; see \cite{Wa08} for details. The
state space of the Kingman $n$-coalescent, $n \in\mathbb{N}$, is the
set $\mathcal{P}_n$ of partitions of the set $\{1,2, \ldots, n\}$. The
process starts in the partition into singletons $\pi_n= \{\{1\}, \ldots,
\{n\}\}$ and has the following dynamics: given that the process is in
the state $\pi_k$, it jumps after a random time $X_k$ to a state $\pi
_{k-1}$ which is obtained by merging two randomly chosen elements from
$\pi_k$. The random inter-coalescence\vspace*{2pt} times $X_k$ are independent,
exponentially distributed random variables with parameters
${k\choose 2}$. The\vspace*{1pt} process can be viewed graphically as a rooted tree
that starts from $n$ leaves labelled from~1 to $n$ and whose any two
branches coalesce independently at rate 1. Each branch of this tree is
situated above a subtree. If this subtree has $r$ leaves, we say that
the branch is \textit{of order} $r$. The branches of order $r\ge2$ are
the \textit{internal} branches, while those of order 1 are the \textit{external} ones (they support subtrees consisting of just one node).

Let us look at the tree from the leaves towards the root
(see Figure \ref{Fig1}). Then the
branch of order $r$ supporting the leaves $i_1, \ldots, i_r$ is formed
at the level ${\sigma}(i_1, \ldots, i_r)$ and ends at level ${\rho
}(i_1, \ldots, i_r)$, where
\[
{\sigma}(i_1, \ldots, i_r) = \max\bigl\{ 1 \leq k \leq
n\dvtx  \{i_1, \ldots, i_r\} \in{\pi}_k\bigr\}
\]
and
\[
{\rho}(i_1, \ldots, i_r) = \max\bigl\{ 1 \leq k < {
\sigma}(i_1, \ldots, i_r)\dvtx  \{i_1, \ldots,
i_r\} \notin{\pi}_k\bigr\}.
\]
For a subset $\{i_1,\ldots,i_r\}$ of leaves, which is not supported by
some branch (which means that $\{i_1,\ldots,i_r\}\notin\pi_k$ for
all $k$) we set ${\sigma}(i_1, \ldots, i_r)={\rho}(i_1, \ldots,\break i_r)=n$.

Let $S_{i_1, \ldots, i_r}$ denote the length of the branch of order $r$
that supports the leaves $i_1, \ldots, i_r$, and write ${\mathcal
L}^{n,r}$ for the total length of order $r$. Then
\[
S_{i_1, \ldots, i_r}= \sum_{l={\rho}(i_1, \ldots, i_r)+1}^{{\sigma
}(i_1, \ldots, i_r)}
X_l
\]
and
\[
{\mathcal L}^{n,r}= \sum_{ 1 \leq i_1< \cdots< i_r \leq n}
S_{i_1, \ldots, i_r}.
\]
Observe that ${\mathcal L}^{n,1}$ is the total length of the external branches.

The length of a randomly chosen external branch in the coalescent tree
has been studied by Freund and M\"ohle \cite{FM09} for the
Bolthausen--Sznitman coalescent and by Gnedin et al. \cite{GIM08} for
the $\Lambda$-coalescent. Asymptotic results concerning the total
external length of Beta$(2-\alpha, \alpha)$-coalescents were given by
M\"ohle \cite{Moe10} for the case $0<\alpha<1$, by Kersting et al.
\cite{KPSJ} for the case $\alpha=1$, and by Kersting et al. \cite
{KSW12} for the case $1<\alpha<2$. For the case $1<\alpha<2$ a weak
law of large numbers result concerning ${\mathcal L}^{n,r}$ can be
easily deduced from Theorem~9 of Berestycki et al. \cite{BBS07} and
also from Dhersin and Yuan \cite{DY12}.

Fu and Li \cite{FL93} computed the expectation and variance of the
total external branch length of the Kingman $n$-coalescent and Caliebe
et al. \cite{CNKR07} derived the asymptotic distribution of a randomly
chosen external branch. In \cite{JK11} Janson and Kersting obtained
the asymptotic normality of the total external branch length. Our main
result states that the same kind of asymptotics holds for the lengths
of order $r \geq1$. Moreover, these lengths are asymptotically independent.

\begin{te*}
For any $s \in\mathbb{N}$, as $n \to\infty$
\[
\sqrt{\frac{n}{4 \log n}} \bigl({\mathcal L}^{n,1}-\mu_1, \ldots, {\mathcal L}^{n,s}-\mu_s \bigr)\stackrel{d} {
\longrightarrow} N (0,I_s ),
\]
where $I_s$ denotes the $s \times s$-identity matrix and $\mu
_r={\mathbb E}({\mathcal L}^{n,r})=\frac{2}{r}$ for every $r \geq1$.
\end{te*}

In a forthcoming paper our theorem will be a main building block for
proving a functional limit theorem for the total external length of the
evolving Kingman-coalescent.

%
%
\begin{figure}

\includegraphics{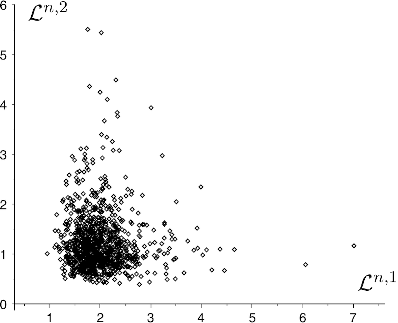}

\caption{External length versus internal length of order 2. The plot
is based on 1000 coalescent realisations with $n=100$.}
\label{Fig2}
\end{figure}

The scatterplot for the lengths of orders 1 and 2 in Figure~\ref{Fig2}
confirms the theorem. The bulk of the points are located around the
mean $(\mu_1,\mu_2)=(2,1)$. Also, in this region hardly any
correlation between the two lengths is visible. The outliers are due to
exceptionally long branches whose occurrence has been explained in
detail in \cite{JK11} for the external case. The simulation shows that
this phenomenon appears similarly in the case of internal lengths, as
one would expect.

As to the proof of the theorem, for the case $s=1$ a hidden symmetry
within the Kingman coalescent is used in \cite{JK11}. Here we
substantially build on the result for $s=1$; however, the proof for the
more general case is rather different. It consists of a coupling device
for Markov chains, which connects the total length of order $r$ to the
total external length: for $1 \leq k \leq n$ let $W_k(r)$ denote the
\textit{number of order $r$ at level} $k$, the number of branches of order
$r$ among the $k$ branches present in the coalescent tree after the
$(n-k)$th coalescing event. (Note that here and elsewhere we are
suppressing the $n$ in the notation.) That is,
\begin{eqnarray*}
W_k(r)&:=& \bigl| \bigl\{ \{i_1, \ldots, i_r\}
\subset\{1, \ldots, n\}\dvtx   i_1< \cdots< i_r,
\\
&&\hspace*{35pt}  {\sigma}(i_1,\ldots, i_r ) \geq k > {
\rho}(i_1,\ldots, i_r) \bigr\} \bigr|.
\end{eqnarray*}
In particular $W_n(r)=0, W_{n-1}(r)=0, \ldots, W_{n-r+2}(r)=0$ and
$W_1(r)=0$ for \mbox{$r<n$}. For an example, see Figure~\ref{Fig1}.

%
%
\begin{figure}

\includegraphics{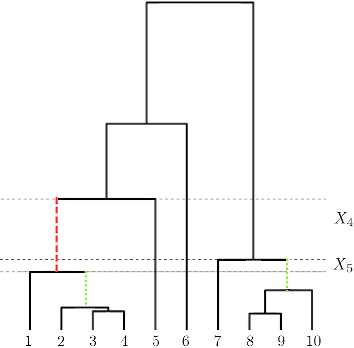}

\caption{The dashed (red) branch is an internal branch of order 4; it
supports the leaves 1, 2, 3, and 4. It is formed at level ${\sigma}(1,
\ldots, 4)=5$ and ends at level ${\rho}(1, \ldots, 4)=3$. Its length
is $S_{1,2,3,4}=X_4+X_5$. The dotted (green) branches are the branches
of order three. The numbers of branches of orders 1 to 10 at level 5
are $W_5(1)=3$, $W_5(2)=0, W_5(3)=1, W_5(4)=1$, and $W_5(i)=0$ for $i
\geq5$.}
\label{Fig1}
\end{figure}

It is important to notice that for any $s \in{\mathbb N}$, the random
vectors $(W_k(1),\break W_k(2),  \ldots, W_k(s))$ form a Markov chain if $k$
runs from $n$ to 1 (a property which facilitated our simulations). The
transition probabilities of the Markov chain are given explicitly in
Section~\ref{sec3}. For a similar approach using a Markov chain embedded in the
Bolthausen--Sznitman coalescent, see \cite{BG08}. The idea of our
proof is to couple $(W_k(r))_{n \geq k \geq1}$ for $1\leq r \leq s$
jointly with $s$ independent copies of the Markov chain of external
numbers $(W_k(1))_{n \geq k \geq1}$. Since in addition the length of
order $r$ is essentially specified by the chain $(W_k(r))_{n \geq k
\geq1}$, it consequently gets the asymptotic behaviour of the external length.

\begin{figure}[t]

\includegraphics{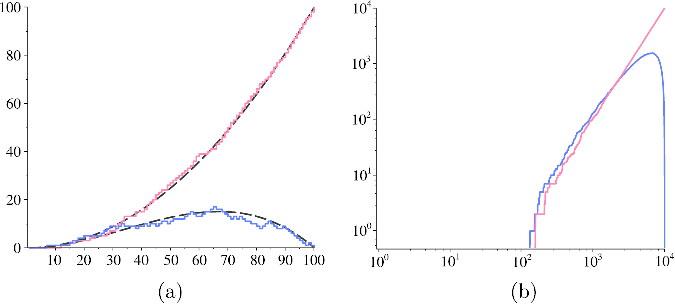}

\caption{\textup{(a)} Simulations of the external numbers $W_k(1)$ (in orange)
and internal numbers $W_k(2)$ of order 2 (in blue) for a coalescent
with $n=100$ for $1 \leq k\leq n$. The black dashed curves represent
the expectations as given in Lemma~\protect\ref{momentsofVandW}. \textup{(b)} Gives the representations in double logarithmic scale for a
coalescent with $n=10^4$.}\label{Fig3}
\end{figure}

The simulations in Figure~\ref{Fig3} give an impression of the
behaviour of the lengths of different orders. In the range between the
levels $n$ and $n^{1-\varepsilon}$ for small ${\varepsilon}>0$ (closer to
the leaves) they differ substantially, as seen in Figure~\ref{Fig3}(a).
This deviation is only due to expectations and does not
appear at the level of fluctuations. Indeed, it is known from \cite
{JK11} that for the external length the fluctuations are induced just
by the $W_k(1)$ with $n^{1-{\varepsilon}} \geq k\geq\sqrt n$. As
suggested by Figure~\ref{Fig3}(b) in this region the evolution of the
chains is similar for orders $r \ge2$. The difference in expectation
is negligible in our construction, as we couple the jumps of the chains
and afterwards consider the lengths of different orders centred at expectation.

The interest in the quantities ${\mathcal L}^{n,r}$ arose from models
where the population is subject to mutation,\vspace*{1pt} the mutations being
modelled as points of a Poisson process with constant rate $\frac
{\theta}{2}$ on\vspace*{1pt} the branches of the coalescent tree. In the infinitely
many sites mutation model, in which every new mutation occurs at a new
locus on the DNA, mutations that are located on the external branches
of the coalescent tree affect only single individuals, whereas
mutations located on an internal branch of order $r\ge2$ affect all
$r$ individuals sitting at the leaves supported by that particular
branch. In a population of size $n$, let $M_r(n)$ denote the number of
mutations carried by exactly $r$ individuals. The vector $(M_1(n), \ldots, M_{n-1}(n))$, called the site frequency spectrum, and the total
number $S_n:= \sum_{r=1}^{n-1} M_r(n)$ of mutations that affect the
population, called the number of segregating sites, are quantities of
statistical importance. Berestycki et al. \cite{BBS07} obtained a weak
law of large numbers for $M_r(n)$, $r \geq1$, in the case of
Beta-coalescents with $1<\alpha<2$.

%
%

For the Kingman coalescent it is known that the number of segregating
sites $S_n$, when rescaled by $\log n$, converges almost surely as $n
\to\infty$ to $\theta$; see, for example, \cite{Be09}, Theorem~2.11.
The expectation of $M_r(n)$ (which is equal to $\frac{\theta}{r}$),
as well as the variances and the covariances of the numbers of
mutations $M_r(n)$, were computed by Fu \cite{Fu95} and Durrett \cite
{Du08}. We obtain the following result as a direct consequence of our theorem.

\begin{cor*}
For any $s \in\mathbb{N}$, as $n \to\infty$
\[
\bigl(M_1(n), \ldots, M_s(n)\bigr) \stackrel{d} {
\longrightarrow} (M_1, \ldots, M_s),
\]
where $M_1, \ldots, M_s$ are independent Poisson-distributed random
variables with parameters $\theta\mu_1, \ldots, \theta\mu_s$.
\end{cor*}

For the proof of the corollary, note from the Poissonian structure of
the mutation process that the characteristic function of $(M_1(n), \ldots, M_s(n))$ is
\begin{eqnarray*}
\varphi_{(M_1(n), \ldots, M_s(n))}(\lambda_1, \ldots, \lambda_s)& =&
{\mathbb E} \bigl[{\mathbb E} \bigl[ e^{i(\lambda_1M_1(n) + \cdots+
\lambda_s M_s(n))} \mid\mathcal{T} \bigr]
\bigr]
\\
& =& {\mathbb E} \bigl[e^{\theta{\mathcal L}^{n,1}(e^{i\lambda_1}-1)}
\cdots e^{\theta{\mathcal L}^{n,s}(e^{i\lambda_s}-1)} \bigr],
\end{eqnarray*}
where $\mathcal{T}$ denotes the ${\sigma}$-algebra containing the
whole information about the coalescent tree. From our theorem it
follows that ${\mathcal L}^{n,r}\stackrel{{\mathbb P}}{\longrightarrow
} \mu_r$ as $n \to\infty$ and therefore
\[
\varphi_{(M_1(n), \ldots, M_s(n))}(\lambda_1, \ldots, \lambda_s)
\longrightarrow e^{\theta\mu_1(e^{i\lambda_1}-1)} \cdots e^{\theta
\mu_s(e^{i\lambda_s}-1)},
\]
as $n \to\infty$.

\begin{rem*}
We note that the convergence ${\mathcal L}^{n,r}\stackrel
{{\mathbb P}}{\longrightarrow} \mu_r$ can also be deduced from the
results of Fu \cite{Fu95}: we have that
\[
{\mathbb V}\bigl(M_r(n)\bigr)={\mathbb V}\bigl({\mathbb E}
\bigl[M_r(n) | \mathcal T\bigr]\bigr)+{\mathbb E}\bigl({\mathbb V}
\bigl[M_r(n) | \mathcal T\bigr]\bigr)={\mathbb V} \biggl(\frac\theta2
{\mathcal L}^{n,r} \biggr)+{\mathbb E} \biggl(\frac\theta2 {\mathcal
L}^{n,r} \biggr).
\]
Comparing this with Fu's formulas (\ref{expectation})--(\ref{alpha2}), we obtain for $r< \frac{n} 2
$ that
\[
{\mathbb V}\bigl({\mathcal L}^{n,r}\bigr)=\frac{2n}{(n-r)(n-r-1)}\sum
_{i=r+1}^n \frac{1} i -
\frac{2} {n-r-1}.
\]
In particular ${\mathbb V}({\mathcal L}^{n,r}) \to0$ and ${\mathcal
L}^{n,r}\stackrel{{\mathbb P}}{\longrightarrow} {\mathbb E}({\mathcal
L}^{n,r})=\mu_r$ as $n \to\infty$.
\end{rem*}

\begin{nota*}
We use the notation $X_n=O_P(f(n))$ for $f(n)>0$ if
\[
\lim_{a \to\infty}\limsup_{n\to\infty}{\mathbb P}
\bigl(X_n > a\cdot f(n)\bigr)=0,
\]
that is, $\frac{X_n}{f(n)}$ is stochastically bounded.\vspace*{2pt}

Throughout $c$ denotes a finite constant whose value is not important
and may change from line to line.
\end{nota*}

\section{Moment computations}\label{sec2}

%
\begin{lm}\label{momentsofVandW}
For the expectation and variance of $W_k(r)$ the following is true. For
$n > r$,
\[
{\mathbb E}\bigl(W_k(r)\bigr) = \frac{(n-k) \cdots(n-k-r+2)}{(n-1) \cdots
(n-r)}\cdot k(k-1)
\quad\mbox{and}\quad
\mathbb{V}\bigl(W_k(r)\bigr)\leq c \frac{k^2}{n},
\]
where $c<\infty$ is a constant depending on $r$. In particular
\[
{\mathbb E}\bigl(W_k(r)\bigr) = \biggl(\frac{n-k}{n}
\biggr)^{r-1}\cdot\frac
{k^2}{n}+O \biggl(\frac{k}{n} \biggr)=
\frac{k^2}{n}+ O \biggl( \frac
{k^3}{n^2}+\frac{k}{n} \biggr) =O
\biggl( \frac{k^2}{n} \biggr)
\]
uniformly in $k \leq n$. Also, for any integer $\alpha\geq2$,
\[
{\mathbb E}\bigl(W_k^{\alpha}(r)\bigr)= \biggl( \biggl(
\frac{n-k}{n} \biggr)^{r-1}\cdot\frac{k^2}{n}
\biggr)^{\alpha}+O \biggl( \biggl(\frac
{k^2}{n} \biggr)^{\alpha-1}+
\frac{k^2}{n} \biggr)=O \biggl( \frac
{k^{2\alpha}}{n^{\alpha}}+\frac{k^2}{n}
\biggr).
\]
\end{lm}

\begin{pf}In order to compute the moments of $W_k(r)$, let us again
label the leaves of the coalescent tree from $1$ to $n$ and note that
$W_{k}(r)$ can be written as
\[
W_k(r)=\sum_{1\leq i_1<\cdots<i_r \leq n} {\mathbf1}_{\{\{i_1, \ldots,
i_r\}\in{\pi}_k\}},
\]
where ${\pi}_k$ is the state of the coalescent process at time $k$.
Then for $n > r$, using the fact that the event $\{\{1, \ldots, r\}\in
{\pi}_k\}$ is the disjoint union (over $n > l_{1} > l_2 > \cdots>
l_{r-1}\geq k$) of the events
$\{$the branch supporting leaves $1, \ldots, r$  is formed
by $r-1$ coalescing events happening at levels $ l_{1}, l_2, \ldots,
l_{r-1}\}$, we have from exchangeability that
%
%
\begin{eqnarray}\label{expectation}
{\mathbb E}\bigl(W_k(r)\bigr)&=& {\mathbb E} \biggl(
\sum_{1\leq i_1<\cdots<i_r
\leq n} {\mathbf1}_{\{\{i_1, \ldots, i_r\}\in{\pi}_k\}} \biggr)
=\pmatrix{n\cr r} {\mathbb P} \bigl(\{1, \ldots, r\}\in{\pi}_k \bigr)\nonumber
\\
&=& \pmatrix{n\cr r} \sum_{n> l_{1} > l_2 > \cdots> l_{r-1}\geq k} \frac
{{n-r\choose 2}}{{n\choose 2}} \cdot
\frac{{n-1-r\choose 2}}{{n-1 \choose 2}} \cdots\nonumber
\\
&&\hspace*{98pt}{}\times
\frac{{l_1+2-r \choose 2} } {{l_1+2\choose 2}} \cdot\frac{{r\choose 2} }{{l_1+1\choose 2}}
\cdot \frac{{l_1-(r-1)\choose  2}}{{l_1\choose 2} } \cdots
\\
&&\hspace*{98pt}{}\times
\frac{{l_j+2-(r-j+1)\choose  2}}
{{l_j+2\choose 2}}
\cdot \frac{{r-j+1\choose  2} }
{{l_j+1\choose 2} } \cdot \frac{ {l_j-(r-j)\choose  2}}
{{l_j\choose  2}} \cdots\nonumber
\\
&&\hspace*{98pt}{}\times \frac{ {l_{r-1}\choose  2}}
{{l_{r-1}+2\choose  2} } \cdot\frac{ {2\choose  2} }{{l_{r-1}+1\choose 2}} \cdot
\frac{{l_{r-1}-1\choose  2} }{{l_{r-1}\choose 2} }
\cdots\frac{ {k\choose  2} }{{k+1\choose  2}}.
\nonumber
\end{eqnarray}
Most binomials in the nominator and the denominator cancel. The
summands turn out to be equal such that
\begin{eqnarray*}
{\mathbb E}\bigl(W_k(r)\bigr)&=&\pmatrix{n\cr r} \sum
_{n> l_{1} > l_2 > \cdots>
l_{r-1}\geq k} \frac{ {r\choose 2} \cdots {2\choose  2} }{{n\choose 2}
\cdots {n-r+1 \choose 2}} \cdot \pmatrix{k\cr 2}
\\
&=&\pmatrix{n\cr r} \pmatrix{n-k\cr r-1} \frac{{r\choose 2} \cdots
{2\choose  2}}{{n\choose  2} \cdots {n-r+1\choose 2} } \cdot \pmatrix{k\cr  2}
\\
&=&\frac{(n-k) \cdots(n-k-r+2)}{(n-1) \cdots(n-r)}\cdot k(k-1).
\nonumber
\end{eqnarray*}
This is the first claim, which directly implies the first asymptotic
formula for ${\mathbb E}(W_k(r))$. Now the second follows by means of
the Bernoulli inequality:
\[
1- \biggl(\frac{n-k} n \biggr)^{r-1}=1- \biggl(1-
\frac{k} n \biggr)^{r-1}\leq(r-1)\frac{k} n.
\]

The computation of the second moment of $W_k(r)$ follows in a similar
way. Note that the event $\{\{i_1, \ldots, i_r\}, \{j_1, \ldots, j_r\}\in
{\pi}_k\}$ is nonempty only if the sets $\{i_1, \ldots, i_r\}$ and $\{
j_1, \ldots, j_r\}$ are identical or disjoint. Thus, for $n > 2r$
\begin{eqnarray*}
{\mathbb E}\bigl(W^2_k(r)\bigr)
&=&\pmatrix{n\cr  r} {\mathbb E} \bigl( {\mathbf1}^2_{\{\{1, \ldots, r\}\in
{\pi}_k\}}
\bigr)
\\
&&{}+ \pmatrix{n\cr  r,r,n-2r} {\mathbb E} ( {\mathbf1}_{\{\{1, \ldots,
r\}\in{\pi}_k\}}\cdot{
\mathbf1}_{\{\{r+1, \ldots,
2r\}\in{\pi}_k\}} )
\\
&=&\pmatrix{n\cr r} {\mathbb P} \bigl( \{1, \ldots, r\}\in{\pi}_k
\bigr)
\\
&&{}+ \pmatrix{n\cr r,r,n-2r} {\mathbb P} \bigl( \{1, \ldots, r\}, \{
r+1, \ldots, 2r\}
\in{\pi}_k \bigr)
\\
&=& {\mathbb E}\bigl(W_k(r)\bigr)
\\
&&{}+ \pmatrix{n\cr  r,r,n-2r}\sum
\frac{{n-2r \choose 2} }{{n\choose  2} } \cdots\frac{ {l''_1+2-2r\choose  2} }{{l''_1+2\choose  2} }
\cdot\frac{{r\choose  2}}{{l''_1+1\choose  2}} \cdot
\frac{{l''_1-(2r-1)\choose  2} }{{l''_1\choose  2}}
\cdots
\\
&&\hspace*{92pt}{}\times \frac{{l''_{2r-2}-1\choose 2}}{{l''_{2r-2}+2\choose 2} } \cdot\frac{1}{{l''_{2r-2}+1\choose 2}}
\cdot\frac{{l''_{2r-2}-2\choose  2}}{{l''_{2r-2}\choose  2} }
\cdots\frac{{k-1\choose  2} }{{k+1\choose  2} },
\end{eqnarray*}
where the sum is taken over all $n> l_{1} > l_2 > \cdots> l_{r-1}\geq
k$ and all $n> l'_{1} > l'_2 > \cdots>l'_{r-1}\geq k$ such that $\{
l_{1}, \ldots, l_{r-1}\} \cap\{l'_{1}, \ldots, l'_{r-1}\} = \varnothing
$. The sequences $(l_j)_{1 \leq j\leq r-1}$ and $(l'_j)_{1 \leq j\leq
r-1}$ denote the coalescence times of the branches supporting leaves
from the sets $\{1, \ldots, r\}$ and $\{r+1, \ldots, 2r\}$, respectively.
The sequence $(l''_j)_{1 \leq j\leq2r-2}$ is the reordering of $l_{1}, \ldots, l_{r-1}, l'_{1}, \ldots, l'_{r-1}$ in decreasing order. Thus
\begin{eqnarray*}
{\mathbb E}\bigl(W^2_k(r)\bigr)&=& {\mathbb E}
\bigl(W_k(r)\bigr)
\\
&&{} +\pmatrix{n\cr  r,r,n-2r} \pmatrix{n-k\cr  r-1, r-1,
n-k-2r+2} \frac{ ({r\choose  2}\cdots
{2\choose  2} )^2}{{n\choose  2}\cdots {n-2r+1\choose  2}}
\\
&&\quad{}\times \pmatrix{k\cr  2} \pmatrix{k-1\cr 2}
\\
&=& {\mathbb E}\bigl(W_k(r)\bigr)
\\
&&{}+ \frac{(n-k)\cdots
(n-k-2r+3)}{(n-1)\cdots
(n-2r)}k(k-1)^2(k-2).
\end{eqnarray*}
The variance of $W_k$ is then for $k \leq n-1$
\begin{eqnarray*}
\mathbb{V}\bigl(W_k(r)\bigr)
&=& {\mathbb E} \bigl(W_k(r) \bigr) \biggl(1+ \frac{(n-k-r+1) \cdots
(n-k-2r+3)}{(n-r-1)\cdots(n-2r)}(k-1)
(k-2)
\\
&&\hspace*{117pt}{}- \frac{(n-k) \cdots
(n-k-r+2)}{(n-1)\cdots(n-r)}\cdot k(k-1) \biggr)
\\
&\leq& {\mathbb E} \bigl(W_k(r) \bigr) \biggl(1+ k(k-1) (n-k)
\cdots(n-k-r+2)
\\
&&\hspace*{64pt}{}\times\biggl( \frac{1}{(n-r-1)\cdots(n-2r)}- \frac{1}{(n-1)\cdots(n-r)}
\biggr) \biggr)
\\
&\leq& {\mathbb E} \bigl(W_k(r) \bigr) \biggl(1+ k^2
n^{r-1} \biggl(\frac{1} {
(n-2r)^r}- \frac{1} {n^r} \biggr) \biggr).
\end{eqnarray*}
Using the mean value theorem we obtain that
\[
\mathbb{V}\bigl(W_k(r)\bigr) \leq{\mathbb E} \bigl(W_k(r)
\bigr) \biggl(1+ k^2 n^{r-1}\frac{2r^2} {(n-2r)^{r+1}} \biggr)\leq c
\frac{k^2}{n},
\]
for $c < \infty$ depending on $r$.

For the other claims we use the same type of argument as above. We have that
%
%
\begin{eqnarray}
\label{alpha1}
&&\pmatrix{n\cr  r, \ldots,r,n-\alpha r} {\mathbb P} \bigl( \{
1, \ldots, r
\}, \ldots, \bigl\{(\alpha-1)r+1, \ldots, \alpha r\bigr\}\in{\pi}_k
\bigr)\nonumber
\\
&&\qquad = \pmatrix{n\cr  r, \ldots,r,n-\alpha r} \pmatrix{n-k\cr  r-1, \ldots, r-1, n-k-
\alpha(r-1)}\nonumber
\\
&&\quad\qquad{}\times \frac{ ({r\choose
2}\cdots{2\choose  2} )^\alpha}{{n\choose  2}
\cdots{n-\alpha r+1\choose  2} } \cdot \pmatrix{k\cr  2}
\cdots \pmatrix{k-\alpha+1\cr 2}\nonumber
\\
&&\qquad = \frac{(n-k)\cdots(n-k-\alpha( r -1)+1)}{(n-1)\cdots(n-\alpha
r)}
\nonumber\\[-8pt]\\[-8pt]
&&\quad\qquad{}\times  k(k-1)^2\cdots(k-\alpha
+1)^2(k-\alpha)
\nonumber
\\
&&\qquad = \frac{(n-k)^{\alpha(r-1)}}{n^{\alpha r}}\cdot k^{2\alpha} + O
\biggl(
\frac{(n-k)^{\alpha(r-1)}}{n^{\alpha r}}\cdot k^{2\alpha-1
} \biggr)
\nonumber
\\
&&\quad\qquad{} + O \biggl( \frac{(n-k)^{\alpha(r-1)-1}}{n^{\alpha
r}}\cdot k^{2\alpha} \biggr)+
O \biggl( \frac{(n-k)^{\alpha
(r-1)}}{n^{\alpha r+1}}\cdot k^{2\alpha} \biggr)
\nonumber
\\
&&\qquad = \frac{(n-k)^{\alpha(r-1)}}{n^{\alpha r}}\cdot k^{2\alpha} + O \biggl
( \frac{k^{2\alpha-1 }}{n^{\alpha}}
\biggr).
\nonumber
\end{eqnarray}
In particular this gives the asymptotic expansion of ${\mathbb
E}(W_k(r))$. Also
%
%
\begin{eqnarray}\label{alpha2}
&& \pmatrix{n\cr  r, \ldots,r,n-\alpha r} {\mathbb P} \bigl( \{1, \ldots, r
\}, \ldots, \bigl\{(\alpha-1)r+1, \ldots, \alpha r\bigr\}\in{\pi}_k
\bigr)
\nonumber\\[-8pt]\\[-8pt]
&&\qquad = O \biggl( \frac{k^{2\alpha}}{n^{\alpha}} \biggr).\nonumber
\end{eqnarray}
Moreover, by expanding $ (\sum_{1\leq i_1<\cdots<i_r \leq n}
{\mathbf1}_{\{\{i_1, \ldots, i_r\}\in{\pi}_k\}} )^\alpha$
%
%
\begin{eqnarray}
\label{alpha3}  {\mathbb E}\bigl(W^\alpha_k(r)\bigr) &=&
\pmatrix{n\cr r, \ldots,r,n-\alpha r} {\mathbb P} \bigl( \{1, \ldots, r\}, \ldots, \bigl
\{(\alpha-1)r+1, \ldots, \alpha r\bigr\}\in{\pi}_k \bigr)\nonumber
\\
&&{} + O \Biggl( \sum_{\beta=1} ^{\alpha-1} \pmatrix{n\cr
r, \ldots,r,n-\beta r}
\\
&&\hspace*{45pt}{}\times  {\mathbb P} \bigl( \{1, \ldots, r\}, \ldots, \bigl\{
(\beta
-1)r+1, \ldots, \beta r\bigr\}\in{\pi}_k \bigr) \Biggr).
\nonumber
\end{eqnarray}
The last claim now follows from (\ref{alpha1}), (\ref{alpha2}),
(\ref{alpha3}), and the fact that $\sum_{\beta=1} ^{\alpha-1}
(\frac{k^{2}}{n} )^{\beta} = O ( (\frac{k^2}{n}
)^{\alpha-1}+\frac{k^2}{n} )$.
\end{pf}

\begin{rem*}
It can be read off from the computation in (\ref
{expectation}) that in the case \mbox{$r=2$}, given the event $\{\rho(1,2)< k
\}=\{\{1,2\} \in\pi_k \}$, the random variable ${\sigma}(1,2)$ is
uniformly distributed on the set of levels $\{k, \ldots, n-1\}$. Indeed
for $k \leq l <n$ it holds
\begin{eqnarray*}
&& {\mathbb P}\bigl(\sigma(1,2) =l, \rho(1,2)< k\bigr)
\\
&&\qquad = \frac{{n-2\choose  2}}{{n\choose  2}} \cdot\frac{
{n-3\choose 2}}{{n-1\choose  2}} \cdots\frac{{l+1\choose  2} } {{l+3\choose  2}} \cdot
\frac{{l\choose  2} } {{l+2\choose  2}} \cdot\frac
{1 }{{l+1\choose  2}} \cdot\frac{{l-1\choose  2}}{{l\choose 2}
}\cdot
\frac{{l-2\choose 2}}{{l-1\choose  2} } \cdots\frac{{k\choose 2} }{ {k+1\choose  2}}
\\
&&\qquad = \frac{1}{{n\choose 2}} \cdot\frac{1}{{n-1\choose 2}}\cdot \pmatrix{k\cr 2},
\end{eqnarray*}
which does not depend on $l$. A similar observation can be made for the
case $r>2$.
\end{rem*}

Using the numbers $W_k(r)$ we have the following simplified expression
for the length of order $r$:
%
%
\begin{equation}
{\mathcal L}^{n,r}=\sum_{2 \leq k \leq n}
W_k(r) \cdot X_k.
\end{equation}

We note from this representation that there are two sources of
randomness in the length of order $r$, one coming from the numbers
$W_k(r)$ and one coming from the exponential inter-coalescence times.
It is easy to see that taking out the randomness introduced by the
exponential times (i.e., replacing them by their expectations) leads to
an error that is asymptotically $O_P (n^{-\nicefrac{1}{2}}
)$ and therefore converges to 0 after the rescaling by $\sqrt{\frac
{n}{\log n}}$. Indeed, by using the independence between the $X_k$'s
and the $W_k(r)$'s and Lemma~\ref{momentsofVandW}, we have that
for a constant $c < \infty$
\begin{eqnarray*}
&& \mathbb{V} \biggl(\sum_{2 \leq k \leq n} W_k(r)
\cdot\bigl( X_k - {\mathbb E}(X_k) \bigr) \biggr)
\\
&&\qquad = \sum
_{2 \leq k\leq n} \mathbb{V} \bigl( W_k(r) \cdot
\bigl( X_k - {\mathbb E}(X_k) \bigr) \bigr)
\\
&&\qquad = \sum_{2 \leq k \leq n} {\mathbb E}\bigl(W_k^2(r)
\bigr) \cdot{\mathbb E}\bigl(\bigl(X_k - {\mathbb E}(X_k)
\bigr)^2\bigr)
\\
&&\qquad  \leq c \sum_{2 \leq k\leq n} \biggl( \frac{k^4}{n^2} +
\frac{k^2}{n} \biggr) \cdot\frac{1}{{k\choose 2^2}} \leq c \frac{1}{n}
\end{eqnarray*}
and therefore
\[
{\mathcal L}^{n,r} =L^{n,r}+ O_P
\bigl(n^{-\nicefrac{1}{2}}\bigr),
\]
where
%
%
\begin{equation}
\label{L} L^{n,r}:= \sum_{2 \leq k \leq n}
W_k(r) \cdot{\mathbb E}(X_k) = \sum
_{2 \leq k \leq n} W_k(r) \cdot\frac{2}{k(k-1)}.
\end{equation}

As a consequence, in the proof of our theorem we need only focus on the
length $L^{n,r}$ which we will, for convenience, still call the length
of order $r$.

\section{The coupling}\label{sec3}

Our proof follows a coupling argument which substantially relies on the
observation that for every $s \in\mathbb{N}$ the vector
$V_k:=(W_k(1), W_k(2), \ldots, W_k(s))$ follows for $n\geq k\geq1$ the
dynamics of an inhomogenous Markov chain with state space $\mathcal
{X}_{n,s}:=\{0,1, \ldots,n\}^s$. We let time run in coalescent direction
(from the leaves to the root of the tree), and for convenience we
consider the evolution of the chain $(V_k)_{n \geq k \geq1}$, running
in the same direction, namely from level $n$ to level 1.

For every $1 \leq r \leq s$ we denote by ${\Delta}W_{n-1}(r), \ldots,
{\Delta}W_{1}(r)$ the sizes of the jumps of the chain $(W_k(r))_{n
\geq k \geq1}$,
\[
{\Delta}W_{k}(r):=W_{k}(r)-W_{k+1}(r),\qquad n-1
\geq k \geq1
\]
and observe that ${\Delta}W_{k}(r) \in\{-2, -1, 0, 1\}$ for all $k$.
The jumps of size 1 correspond to the levels at which a new branch of
order $r$ is formed (by the coalescence of two other branches), whereas
the jumps of sizes $-$1 and $-$2 happen at the levels at which one (or,
resp., two) branches of order $r$ end (by coalescence of one of them
with some other branch or by mutual coalescence).

For $1 \leq k \leq n$ and $v, v' \in\mathcal{X}_{n,s}$ let
\[
P_{v}^k\bigl(v'\bigr):={\mathbb P}
\bigl(V_{k-1}=v' \mid V_{k}=v \bigr)
\]
denote the transition probabilities of $(V_k)_{n \geq k \geq1}$. They
are given for $v=(w_1, \ldots, w_s)$, $w_1+ \cdots+ w_s \leq k$ by
%
%
\begin{equation}
\label{transitions}
 P_{v}^k\bigl(v'\bigr)=
\cases{
\displaystyle \frac{{ {k-w_1-\cdots- w_s\choose 2}} }{{ k\choose 2}},
\cr
\hspace*{32pt}
\qquad \mbox{if $v=v'$,}
\vspace*{5pt}\cr
\displaystyle \frac{{ w_i (k-\sum_{i=1}^s w_i )}}{{ k\choose 2}},
\cr
\hspace*{32pt}\qquad \mbox{if $v=v'-e_{i}$ for some $i$,}
\vspace*{5pt}\cr
\displaystyle\frac{{ {w_i\choose 2}} }{{ k\choose 2}},
\hspace*{7.5pt}\qquad \mbox{if $v=v'-2e_{i}+e_{2i}$ for some $i$,}
\vspace*{5pt}\cr
\displaystyle\frac{{ w_i w_j}}{{ k\choose 2}},
\qquad \mbox{if $v=v'-e_{i}-e_{j}+e_{i+j}$ for some $i\neq j$,}
\vspace*{5pt}\cr
0,\hspace*{21.5pt}
\qquad \mbox{else,}}
\end{equation}
where $e_i=(\delta_{i,l})_{1 \leq l\leq s}$ [note that $e_i=(0, \ldots,
0)$ for $i >s$]. For all other $v \in\mathcal{X}_{n,s}$ we set for
definiteness $P_{v}^k(v')=\delta_{v, v'}$ in order to obtain a proper
transition matrix.

In particular for $s=2$ the transition probabilities in (\ref
{transitions}) are given for $w_1+w_2 \leq k$ by
\begin{eqnarray*}
P_{(w_1,w_2)}^k(w_1,w_2)&=&
\frac{{k-w_1-w_2\choose 2} }{{ k\choose 2}},
\\
P_{(w_1,w_2)}^k(w_1-1,w_2-1)&=& \frac{w_1w_2}{{k\choose 2}},
\\
P_{(w_1,w_2)}^k(w_1-1,w_2)&=&
\frac{w_1(k-w_1-w_2)}{{k\choose 2}},
\\
P_{(w_1,w_2)}^k(w_1-2,w_2+1)&=&\frac{{w_1\choose 2} }{{k\choose 2}},
\\
P_{(w_1,w_2)}^k(w_1,w_2-1)&=& \frac{w_2(k-w_1-w_2)}{{k\choose 2}},
\\
P_{(w_1,w_2)}^k(w_1,w_2-2)&=& \frac{{w_2\choose  2} }{{k\choose 2}}
\end{eqnarray*}
and for $s=1$ for $w \leq k$ by
%
%
\begin{eqnarray}\label{transitions1}
P_w^k(w)&=&\frac{{k-w\choose 2} }{{k\choose 2}},\qquad
P_w^k(w-1)=\frac
{w(k-w)}{{k\choose 2}},
\nonumber\\[-8pt]\\[-8pt]
P_w^k(w-2)&=&\frac{{w\choose  2} }{{k\choose 2}}.\nonumber
\end{eqnarray}

Let us now describe the coupling in detail. Let $1 <a_n \leq n$ and $s
\in{\mathbb N}$ be fixed. Starting\vspace*{2pt} at $a_n$ we couple the Markov chain
$ (V_k )_{a_n\geq k\geq1}$, $V_k=({W}_k(1), \ldots, W_k(s))$,
with another chain $ (\widetilde V _k )_{a_n\geq k\geq1}$,
$\widetilde{V}_k =(\widetilde{W}_k(1), \ldots, \widetilde W_k(s))$,
defined on the same probability space as $ (V_k )_{a_n\geq
k\geq1}$ and having the same state space $\mathcal{X}_{n,s}$. The
components of $ (\widetilde{V}_k )_{a_n\geq k\geq1}$ evolve
as independent copies of $ (W_k(1) )_{a_n\geq k\geq1}$, the
Markov\vspace*{2pt} chain of external numbers. Therefore its transition
probabilities $\widetilde{P}_{v}^k(\cdot)$ for $v \in\mathcal
{X}_{n,s}$ and $1 < k \leq a_n$ are given by the product of the
transition probabilities of its $s$ components, given in (\ref
{transitions1}). The process $ ( (V_k, \widetilde{V}_k )
)_{a_n\geq k\geq1}$ is constructed as a Markov chain, where the
jumps are coupled in a way that we will describe in detail shortly.\vspace*{2pt}

Thus let $Q_{v}^k$ and $\widetilde Q_{\widetilde v}^k$ denote the
conditional distributions of the jumps ${\Delta}{V}_k$ and ${\Delta
}\widetilde{V}_k$ of the two Markov chains, given the current states
$v$ and $\widetilde v$, respectively. (The notation ${\Delta}{V}_k$ and
${\Delta}\widetilde{V}_k$ refer to component-wise differences.) For
the sequel it is important that the leading terms of $Q_{v}^k$ and
$\widetilde Q_{\widetilde v}^k$ agree. More precisely, from~(\ref
{transitions}) and~(\ref{transitions1}), under the constrain that
$w_1+ \cdots+ w_s \leq k$
%
\begin{equation}
\label{transitionsshort}
Q_{v}^k(z), \widetilde Q_{v}^k(z)=
\cases{
\displaystyle 1- {\sum_{j=1}^s} \frac{{ 2w_j}}{{k}} + O \Biggl( {\sum_{j=1}^s}\frac{{w_j^2}}{{ k^2}} \Biggr), & \quad if $z=(0, \ldots, 0)$,
\vspace*{5pt}\cr
\displaystyle \frac{{ 2w_i}}{ {k} }+ O \Biggl( {\sum_{j=1}^s}\frac{{ w_j^2}}{{ k^2}} \Biggr), &\quad if $z=-e_i$ for some $i$,
\vspace*{5pt}\cr
\displaystyle O \Biggl( {\sum_{j=1}^s}\frac{{w_j^2}}{{ k^2}} \Biggr), &\quad else,}\hspace*{-20pt}
\end{equation}
where $z \in\{-2, -1, 0, 1\}^s$ and $e_i=(\delta_{i,l})_{1 \leq l\leq
s}$. Here we use that $w_iw_j \leq w_i^2 +w_j^2$ and $w_i \leq k$.

As it is well known (see, e.g., \cite{LPW09}), an optimal coupling of
the two distributions is specified as follows. Let $\| \cdot\|_{\mathrm{TV}}$
denote the total variation distance between two distributions, and define
\[
p=p_{v, \widetilde v}:=1- \bigl\|Q_{v}^k - \widetilde
Q_{\widetilde v}^k \bigr\|_{\mathrm{TV}}.
\]

Then, with probability $p$ choose ${\Delta}{V}_k={\Delta}\widetilde
{V}_k=Z$, where the random variable $Z$ has distribution $\gamma_{I}$,
given by its weights
\[
\gamma_I(z) = \frac{Q_{v}^k (z) \wedge\widetilde Q_{\widetilde
v}^k(z)}{p},
\]
$z\in\{-2, -1, 0, 1\}^s$, and with probability $1-p$ choose ${\Delta
}{V}_k$ according to the probability distribution weights
\[
\gamma_{\mathit{II}}(z) = \frac{ (Q_{v}^k (z) - \widetilde
Q_{\widetilde v}(z) )^+}{1-p},
\]
and independently choose ${\Delta}\widetilde{V}_k$ according to the
probability distribution weights
\[
\gamma_{\mathit{III}}(z) =\frac{ (\widetilde Q_{\widetilde v}(z)-Q_{v}^k
(z) )^+}{1-p},
\]
$z \in\{-2, -1, 0, 1\}^s$.

This coupling is optimal in the sense that the probability
${\mathbb P} ({\Delta}{V}_k \neq{\Delta}\widetilde{V}_k \mid\break
V_k=v, \widetilde{V}_k=\widetilde v )$
is minimal among the corresponding\vspace*{1pt} probabilities for couplings of the
two distributions $Q_{v}^k$ and $\widetilde Q_{\widetilde v}^k$, and
therefore it is equal to $\|Q_{v}^k - \widetilde Q_{\widetilde v}^k \|
_{\mathrm{TV}}$. As starting distribution of the coupled chain $ (V_k,
\widetilde{V}_k )_{a_n\geq k\geq1}$ we allow any distribution of
$ (V_{a_n}, \widetilde{V}_{a_n} )$ such that the marginals
are the distributions of $V_{a_n}$ and $\widetilde{V}_{a_n}$,
respectively. We point out that the distributions of $V_{a_n}$ and
$\widetilde{V}_{a_n}$ are given by the Kingman coalescent at level
$a_n$. Up to this constraint the common distribution is arbitrary.

The next two lemmas give essential properties of the coupling.
%

\begin{lm}\label{lemmapk}
There is a $c <\infty$ such that the above defined coupling satisfies
for $r \leq s$
\[
{\mathbb P} \bigl( {\Delta}W_{k}(r) \neq{\Delta}\widetilde
{W}_{k}(r) \bigr) \leq c \biggl(\frac{k}{a_n \sqrt{n}}+\frac{a_n
k}{n^2}+
\frac{1} {k} \biggr)
\]
and
\[
{\mathbb E} \bigl( \bigl| W_{k}(r) -\widetilde{W}_{k}(r) \bigr|
\bigr) \leq c \biggl(\frac{k^2}{a_n \sqrt{n}}+\frac{a_n k^2}{n^2}+1
\biggr)
\]
for all $1 \leq k < a_n$.
\end{lm}

\begin{pf}In the proof we write as an abbreviation $W_k$ instead of
$W_k(r)$ and similarly $\widetilde W_k$, ${\Delta}W_k$ and ${\Delta
}\widetilde W_k$ instead of $\widetilde W_k(r)$, ${\Delta}W_k(r)$ and
${\Delta}\widetilde W_k(r)$, respectively.\vspace*{1pt}

From (\ref{transitionsshort}) it follows that for both the chains
$ (V_k )_{a_n\geq k\geq1}$ and $ (\widetilde{V}_k
)_{a_n\geq k\geq1}$ jumps of sizes $(0, \ldots, 0)$ and $-e_i$ with
$1\leq i\leq r$ occur with probabilities of larger order than jumps of
other sizes. It follows from the coupling that
\begin{eqnarray*}
\{{\Delta}W_k \neq{\Delta}\widetilde{W}_{k}\} & \subset&
\{{\Delta}W_k=-1, {\Delta}V_k \neq{\Delta}
\widetilde{V}_k \}
\\
&&{} \cup\{{\Delta}\widetilde{W}_{k}=-1, {
\Delta}V_k \neq{\Delta}\widetilde{V}_k \}
\\
&&{}\cup\bigl\{{\Delta}W_k \in\{1, -2\}\bigr\} \cup\{{
\Delta}\widetilde{W}_{k}=-2\}.
\end{eqnarray*}
Note that since $\widetilde{W}_{k}$ has the distribution of the
external number at level $k$, the jumps size $\Delta\widetilde
{W}_{k}$ cannot take the value 1.

Thus, writing as an abbreviation ${\mathbb P}^k_{v,\widetilde
{v}}(\cdot)$ for the conditional probability given the event $\{V_k=v,
\widetilde{V}_k = \widetilde{v}\}$, we obtain
\begin{eqnarray*}
{\mathbb P}^{k+1}_{v,\widetilde{v}} ( {\Delta}W_{k} \neq{
\Delta}\widetilde{W}_{k} )& \leq& (1-p) \gamma_{\mathit{II}}({
\Delta}W_k=-1) + (1-p) \gamma_{\mathit{III}}({\Delta}\widetilde
W_k=-1)
\\
&&{} + {\mathbb P}^{k+1}_{v,\widetilde{v}} \bigl( {\Delta
}W_{k}\in\{1, -2\} \bigr) + {\mathbb P}^{k+1}_{v,\widetilde{v}}
( {\Delta}\widetilde{W}_{k} =-2 )
\nonumber
\\
& \leq& (1-p) \gamma_{\mathit{II}}( {\Delta}V_k=-e_r) +
(1-p) \gamma_{\mathit{III}}({\Delta}\widetilde{V}_k=-e_r)
\nonumber
\\
&&{} + c \cdot\sum_{i=1}^s
\frac{ w_i^2+\widetilde
{w}_i^2}{k^2}
\nonumber
\\
& &{}+ {\mathbb P}^{k+1}_{v,\widetilde{v}} \bigl( {
\Delta}W_{k}\in\{1, -2\} \bigr) + {\mathbb P}^{k+1}_{v,\widetilde
{v}}
( {\Delta}\widetilde{W}_{k} =-2 )
\nonumber
\\
& \leq& (1-p) \gamma_{\mathit{II}}( {\Delta}V_k=-e_r) +
(1-p) \gamma_{\mathit{III}}({\Delta}\widetilde{V}_k=-e_r)
\nonumber
\\
&&{} + c \cdot\sum_{i=1}^s
\frac{ w_i^2+\widetilde
{w}_i^2}{k^2}.
\nonumber
\end{eqnarray*}
Using now the definitions of $\gamma_{\mathit{II}}$ and $\gamma_{\mathit{III}}$ we get that
%
\begin{equation}
\label{probneq2} \qquad {\mathbb P}^{k+1}_{v,\widetilde{v}} ( {
\Delta}W_{k} \neq{\Delta}\widetilde{W}_{k} ) \leq\bigl|
Q_{v}^{k+1}(-e_r) - \widetilde
Q_{\widetilde{v}}^{k+1}(-e_r) \bigr| + c \cdot\sum
_{i=1}^s \frac{
w_i^2+\widetilde{w}_i^2}{k^2}.
\end{equation}
Let us introduce the filtration $\mathbb{F}=(\mathcal{F}_k)_{1\leq
k\leq a_n}$ with $\mathcal{F}_{a_n} \subset\mathcal{F}_{a_n-1}
\subset\cdots\subset\mathcal{F}_1$ defined by
\[
\mathcal{F}_k = \sigma\bigl((V_j)_{k\leq j \leq a_n}, (
\widetilde{V}_j)_{k\leq j \leq a_n} \bigr).
\]
Then (\ref{probneq2}) in view of (\ref{transitionsshort}) may be
written as
%
\begin{eqnarray}
\label{pkcond}
&& {\mathbb P} ({\Delta}W_{k} \neq{\Delta}\widetilde
{W}_{k} \mid{\mathcal F}_{k+1} )
\nonumber\\[-8pt]\\[-8pt]
&&\qquad  \leq\frac{2}{k}|W_{k+1}-\widetilde W_{k+1}| + c \sum
_{i=1}^s \frac
{ W_{k+1}^2(i)+\widetilde W_{k+1}^2(i)}{k^2}
\nonumber
\end{eqnarray}
for a constant $c < \infty$. Taking expectation in the inequality
above, we obtain using Lemma~\ref{momentsofVandW}
%
%
\begin{equation}
\label{pk} {\mathbb P} ({\Delta}W_{k} \neq{\Delta}
\widetilde{W}_{k} ) \leq\frac{2}{k}\cdot{\mathbb E}
\bigl(|W_{k+1}-\widetilde W_{k+1}| \bigr) + c \biggl(\frac{k^2}{n^2}
+\frac{1}{n} \biggr).
\end{equation}

We now proceed to finding a bound for ${\mathbb E}
(|W_{k}-\widetilde W_{k}| )$ for $2 \leq k \leq a_n$. From the
transition probabilities (\ref{transitions}) we get that
\begin{eqnarray*}
{\mathbb E}[{\Delta}W _{k} \mid\mathcal{F}_{k+1} ]
&=&(-1) \cdot\frac{W_{k+1}(r)(k+1-W_{k+1}(1)- \cdots-
W_{k+1}(r))}{{k+1\choose  2} }
\nonumber
\\
&&{}+ (-1) \cdot\frac{W_{k+1}(r)W_{k+1}(1)+ \cdots
+W_{k+1}(r) W_{k+1}(r-1) }{{k+1\choose  2} }
\nonumber
\\
&&{} + (-2) \cdot\frac{{W_{k+1}(r)\choose 2} }{{k+1\choose 2} } + 1 \cdot\frac{Z_{k+1}}{{k+1\choose 2} }
\nonumber
\\
&=& -\frac{2}{k+1}W_{k+1} + \frac{Z_{k+1}}{{k+1\choose 2}},
\nonumber
\end{eqnarray*}
where, letting $d_r=1$ if $r$ is even and 0 otherwise,
%
%
\begin{equation}
\label{defZ} Z_k=Z_{k}(r):= \mathop{\sum_{1 \leq i\leq r-1}}_{i \neq r-i} W_{k}(i)W_{k}(r-i)+
d_r \cdot\pmatrix{W_{k}(\nicefrac{r} {2})\cr 2}.
\end{equation}
Therefore
%
%
\begin{equation}
\label{deltaWk} {\mathbb E}[{\Delta}W_{k} \mid\mathcal{F}_{k+1}
] = -\frac
{2}{k+1}W_{k+1} + \frac{Z_{k+1}}{{k+1\choose 2}}
\end{equation}
and also with a similar but even simpler calculation using (\ref
{transitions1}),
%
%
\begin{equation}
\label{deltaWktilde} {\mathbb E}[{\Delta}\widetilde W_{k} \mid{\mathcal
F}_{k+1}]=-\frac
{2}{k+1} \widetilde W_{k+1}.
\end{equation}
Now note that the absolute value of the difference between the jumps of
$W_{k+1}$ and $\widetilde W_{k+1}$ is at most 3. Thus
\begin{eqnarray*}
&& {\mathbb E} \bigl[|W_{k}-\widetilde W_{k} | \mid{\mathcal
F}_{k+1}\bigr]
\\
&&\qquad ={\mathbb E}\bigl[|W_{k+1}-\widetilde W_{k+1} + {
\Delta}W_k -{\Delta}\widetilde W_{k} | \mid{\mathcal
F}_{k+1}\bigr]
\\
&&\qquad  \leq{\mathbb E}[W_{k+1}-\widetilde W_{k+1} + {
\Delta}W_k -{\Delta}\widetilde W_{k} \mid{\mathcal
F}_{k+1}] \cdot{\mathbf1}_{\{
W_{k+1}-\widetilde W_{k+1} \geq3\}}
\\
&&\quad\qquad{} +{\mathbb E}[\widetilde W_{k+1}-W_{k+1} + {
\Delta}\widetilde W_k -{\Delta}W_{k} \mid{\mathcal
F}_{k+1}] \cdot{\mathbf1}_{\{W_{k+1}-\widetilde W_{k+1} \leq-3\}}
\\
&&\quad\qquad{} + \bigl(|W_{k+1}-\widetilde W_{k+1}| + {
\mathbb E}\bigl[| {\Delta}W_k -{\Delta}\widetilde W_{k} | \mid{
\mathcal F}_{k+1}\bigr] \bigr) \cdot{\mathbf1}_{\{|W_{k+1}-\widetilde
W_{k+1}| \leq2\}}.
\end{eqnarray*}
Using (\ref{deltaWk}) and (\ref{deltaWktilde}) we obtain
\begin{eqnarray*}
&& {\mathbb E}\bigl[|W_{k}-\widetilde W_{k} | \mid{\mathcal
F}_{k+1}\bigr]
\\
&&\qquad  \leq\biggl(W_{k+1}-\widetilde W_{k+1} - \frac{2}{k+1}
(W_{k+1}-\widetilde W_{k+1} )+ \frac{Z_{k+1}}{{k+1\choose 2}} \biggr)
\cdot{\mathbf1}_{\{W_{k+1}-\widetilde W_{k+1} \geq3\}}
\\
&&\quad\qquad{} + \biggl(\widetilde W_{k+1}-W_{k+1} -
\frac{2}{k+1} (\widetilde W_{k+1}-W_{k+1} ) -
\frac{Z_{k+1}}{{ k+1\choose
2}} \biggr) \cdot{\mathbf1}_{\{W_{k+1}-\widetilde W_{k+1} \leq-3\}}
\\
&&\quad\qquad{} + \bigl(|W_{k+1}-\widetilde W_{k+1}| + 3 \cdot{
\mathbb P}( {\Delta}W_k \neq{\Delta}\widetilde W_{k} \mid{
\mathcal F}_{k+1}) \bigr) \cdot{\mathbf1}_{\{|W_{k+1}-\widetilde
W_{k+1}| \leq
2\}}.
\end{eqnarray*}
By (\ref{pkcond}) we have that
\begin{eqnarray*}
&& {\mathbb E}\bigl[|W_{k}-\widetilde W_{k} | \mid{\mathcal
F}_{k+1}\bigr]
\\
&&\qquad \leq\biggl(|W_{k+1}-\widetilde W_{k+1} |-
\frac
{2}{k+1}|W_{k+1}-\widetilde W_{k+1} | +
\frac{Z_{k+1}}{{k+1\choose 2}} \biggr) \cdot{\mathbf1}_{\{
|W_{k+1}-\widetilde W_{k+1}| \geq3\}}
\\
&&\quad\qquad{}  + \Biggl(|W_{k+1}-\widetilde W_{k+1} |+
\frac
{6}{k}|W_{k+1}-\widetilde W_{k+1} |
\\
&&\hspace*{53pt}\quad\qquad{} +c \sum_{i=1}^s
\frac{
W_{k+1}^2(i)+\widetilde W_{k+1}^2(i)}{k^2} \Biggr) \cdot{\mathbf1}_{\{
|W_{k+1}-\widetilde W_{k+1}| \leq2\}}
\\
&&\qquad \leq|W_{k+1}-\widetilde W_{k+1}| \biggl(1-\frac{2}{k+1}
\biggr) + \frac{16}{k}
\\
&&\quad\qquad{}  + c \sum_{i=1}^s
\frac{
W_{k+1}^2(i)+\widetilde W_{k+1}^2(i)}{k^2} + \frac{Z_{k+1}}{{k+1\choose 2}}.
\end{eqnarray*}

Taking expectation and using Lemma~\ref{momentsofVandW} and the
fact that $k \leq n$, we obtain that
\[
{\mathbb E}\bigl(|W_{k}-\widetilde W_{k} |\bigr) \leq\biggl(1-
\frac{2}{k+1} \biggr){\mathbb E}\bigl[|W_{k+1}-\widetilde
W_{k+1} |\bigr]+ c \biggl(\frac{k^2} {n^2} +\frac{1} {k} \biggr).
\nonumber
\]

Dividing the previous inequality by $k(k-1)$ we obtain a recurrence
formula that we iterate from $k$ up to $a_n-1$,
\begin{eqnarray*}
&& \frac{1}{k(k-1)}  {\mathbb E}\bigl(|W_{k} -\widetilde
W_{k}|\bigr)
\\
&&\qquad  \leq \frac{1}{k(k+1)} \biggl( {\mathbb E}\bigl(|W_{k+1}-\widetilde
W_{k+1}|\bigr)+ c \biggl( \frac{k^2} {n^2} + \frac{1} {k} \biggr)
\biggr)
\\
&&\qquad  \leq \frac{1}{a_n(a_n-1)} {\mathbb E}\bigl(|W_{a_n}-\widetilde
W_{a_n}|\bigr)+ c \sum_{j=k}^{a_n-1}
\biggl( \frac{1} {n^2} + \frac{1} {j^3} \biggr)
\\
&&\qquad  \leq \frac{1}{a_n(a_n-1)} \bigl({\mathbb E} \bigl(\bigl|W_{a_n}-{\mathbb
E}(W_{a_n})\bigr| \bigr)+ \bigl|{\mathbb E}(W_{a_n})-{\mathbb E}(
\widetilde W_{a_n}) \bigr|
\\
&&\hspace*{173pt}{}+{\mathbb E} \bigl(\bigl|\widetilde W_{a_n} - {
\mathbb E}(\widetilde W_{a_n})\bigr| \bigr) \bigr) + c \biggl(
\frac{a_n} {
n^2} + \frac{1} {k^2} \biggr).
\end{eqnarray*}
Finally by Lemma~\ref{momentsofVandW},
\begin{eqnarray*}
\frac{1}{k(k-1)} {\mathbb E}\bigl(|W_{k} -\widetilde W_{k}|\bigr)
&\leq&  c \biggl( \frac{1}{a_n(a_n-1)} \biggl(\frac{a_n}{\sqrt{n}}+
\frac
{a_n^{3}}{n^2}+\frac{a_n}{n} \biggr)+\frac{ a_n} {n^2} +
\frac{1} {
k^2} \biggr)
\\
& \leq& c \biggl( \frac{1}{a_n \sqrt{n}}+\frac{a_n}{n^2}+\frac{1} {
k^2}
\biggr).
\end{eqnarray*}
This gives the second claim of the lemma. Using this claim and the fact
that $1\leq k\leq a_n\leq n$ in (\ref{pk}) yields the first claim.
\end{pf}

%
\begin{lm}\label{lemmavar}
There is a constant $c<\infty$ such that for $r \leq s$ it holds that
\[
\mathbb{V} \bigl( W_{k}(r) - \widetilde W_{k}(r) \bigr)
\leq c \cdot\biggl(\frac{k^2}{a_n\sqrt{n} } + \frac{a_n k^2} {n^2} +
\frac
{k^3}{a_n n}+1 \biggr)
\]
for all $1 \leq k < a_n$.
\end{lm}

\begin{pf}
We\vspace*{2pt} again write here as an abbreviation $W_k$, $\widetilde
W_k$, ${\Delta}W_k$, and ${\Delta}\widetilde W_k$ instead of
$W_k(r)$, $\widetilde W_k(r)$, ${\Delta}W_k(r)$, and ${\Delta
}\widetilde W_k(r)$, respectively.

Using (\ref{deltaWk}) and (\ref{deltaWktilde}) together with the
fact that $|{\Delta}W_{k-1} - {\Delta}\widetilde W_{k-1}|\leq3$, we obtain
\begin{eqnarray*}
&& \mathbb{V} ( W_{k-1} - \widetilde W_{k-1} )
\\
&&\qquad  = \mathbb{V} (
W_{k} - \widetilde W_{k} + {\Delta}W_{k-1} - {
\Delta}\widetilde W_{k-1} )
\\
&&\qquad  \leq \mathbb{V} ( W_{k} - \widetilde W_{k} )
\\
&&\quad\qquad {}
+ 2 {\mathbb E} \bigl({\mathbb E} \bigl[ \bigl(W_{k}-
\widetilde W_{k} - {\mathbb E} ( W_{k}- \widetilde
W_{k} ) \bigr)
\\
&&\hspace*{71pt}{}\times \bigl({\Delta}W_{k-1} - {\Delta}
\widetilde W_{k-1}-{\mathbb E} ( {\Delta}W_{k-1} - {\Delta}
\widetilde W_{k-1} ) \bigr) | {\mathcal F}_k \bigr] \bigr)
\\
&&\quad\qquad{} +{\mathbb E} ({\Delta}W_{k-1} - {\Delta}\widetilde
W_{k-1} )^2
\nonumber
\\
&&\qquad = \biggl(1-\frac{4}{k} \biggr) \mathbb{V} ( W_{k} -
\widetilde W_{k} )
\\
&&\quad\qquad{} + 2 {\mathbb E} \biggl( \bigl(W_{k}- \widetilde
W_{k}- {\mathbb E} (W_{k}- \widetilde W_{k} )
\bigr) \cdot\frac
{Z_k-{\mathbb E}(Z_k)}{{k\choose 2}} \biggr)
\\
&&\quad\qquad{}+ 9 {\mathbb P}(W_{k-1} \neq{\Delta}\widetilde
W_{k-1}).
\nonumber
\end{eqnarray*}
Applying now the Cauchy--Schwarz inequality for the second term and
Lemma~\ref{lemmapk} for the third term on the left-hand side of the
inequality above, we obtain that for a constant $c < \infty$
%
\begin{eqnarray}
\label{variances}
&& \mathbb{V} ( W_{k-1} - \widetilde W_{k-1} )\nonumber
\\
&&\qquad \leq\biggl(1-\frac{4}{k} \biggr) \mathbb{V} ( W_{k} -
\widetilde W_{k} ) + \frac{4}{k(k-1)} \bigl(\mathbb{V}
(W_{k}- \widetilde W_{k} ) \bigr)^{\nicefrac{1}{2}}\cdot
\bigl(\mathbb{V} (Z_k) \bigr)^{\nicefrac{1}{2}}
\\
&&\quad\qquad{} + c \biggl(\frac{k}{a_n \sqrt{n}}+\frac{a_n k}{n^2}+\frac{1} {
k} \biggr).\nonumber
\end{eqnarray}
Let us now look closer at the variance of $Z_k$. In order to bound it
from above, it is sufficient to bound the terms of the form $\mathbb
{V}(W_k(i)W_k(j))$, $1 \leq i,j\leq r$; see the definition of $Z_k$ in
(\ref{defZ}). Writing as an abbreviation $W'_k$ and $W''_k$ for
$W_k(i)$ and $W_k(j)$, respectively, we have that
\begin{eqnarray*}
\mathbb{V}\bigl(W_k(i)W_k(j)\bigr) &\leq& {\mathbb E}
\bigl(W'_kW''_k
- {\mathbb E}\bigl(W'_k\bigr){\mathbb E}
\bigl(W''_k\bigr) \bigr)^2
\\
&=& {\mathbb E} \bigl(\bigl(W'_k - {\mathbb E}
\bigl(W'_k\bigr)\bigr) W''_k
+ {\mathbb E}\bigl(W'_k\bigr) \bigl(W''_k
- {\mathbb E}\bigl(W''_k\bigr)\bigr)
\bigr)^2
\\
& \leq& 2{\mathbb E} \bigl(\bigl(W'_k - {\mathbb E}
\bigl(W'_k\bigr)\bigr)^2
\bigl(W''_k\bigr)^2 + \bigl({
\mathbb E}\bigl(W'_k\bigr)\bigr)^2
\bigl(W''_k - {\mathbb E}
\bigl(W''_k\bigr)\bigr)^2
\bigr).
\end{eqnarray*}
Using the fact that $W''_k \leq k$ and then Lemma~\ref
{momentsofVandW}, we obtain
%
\begin{eqnarray}\label{WW}
&& \mathbb{V}\bigl(W_k(i)W_k(j)\bigr)\nonumber
\\
&&\qquad \leq 2{\mathbb E} \bigl(\bigl(W'_k - {\mathbb E}
\bigl(W'_k\bigr)\bigr)^2
\bigl(W_k''-{\mathbb E}
\bigl(W''_k\bigr)\bigr)W''_k
\bigr)\nonumber
\\
&&\quad\qquad{} + 2{\mathbb E} \bigl(\bigl(W'_k - {
\mathbb E}\bigl(W'_k\bigr)\bigr)^2 {\mathbb
E}\bigl(W''_k\bigr)W''_k
\bigr) +2\bigl({\mathbb E}\bigl(W'_k\bigr)
\bigr)^2 \mathbb{V}\bigl(W''_k
\bigr)\nonumber
\\
&&\qquad \leq 2{\mathbb E} \bigl(\bigl(W'_k - {\mathbb E}
\bigl(W'_k\bigr)\bigr)^2
\bigl(W''_k-{\mathbb E}
\bigl(W''_k\bigr)\bigr)^2 \bigr)\nonumber
\\
&&\quad\qquad{} + 2{\mathbb E} \bigl(\bigl(W'_k - {
\mathbb E}\bigl(W'_k\bigr)\bigr)^2
\bigl(W''_k-{\mathbb E}
\bigl(W''_k\bigr)\bigr){\mathbb E}
\bigl(W''_k\bigr) \bigr)
\\
&&\quad\qquad{} + 2k \cdot{\mathbb E} \bigl(\bigl(W'_k
- {\mathbb E}\bigl(W'_k\bigr)\bigr)^2 {
\mathbb E}\bigl(W''_k\bigr) \bigr) +2
\bigl( {\mathbb E}\bigl(W_k'\bigr)\bigr)^2
\mathbb{V}\bigl(W''_k\bigr)
\nonumber
\\
&&\qquad \leq 2{\mathbb E} \bigl(\bigl(W'_k - {\mathbb E}
\bigl(W'_k\bigr)\bigr)^2
\bigl(W''_k-{\mathbb E}
\bigl(W''_k\bigr)\bigr)^2 \bigr)\nonumber
\\
&&\quad\qquad{} +4 k \cdot{\mathbb E} \bigl(\bigl(W'_k
- {\mathbb E}\bigl(W'_k\bigr)\bigr)^2 {
\mathbb E}\bigl(W''_k\bigr) \bigr) + 2
\bigl({\mathbb E}\bigl(W'_k\bigr)\bigr)^2
\mathbb{V}\bigl(W''_k\bigr)
\nonumber
\\
&&\qquad \leq2 {\mathbb E}  \bigl(\bigl(W'_k - {\mathbb E}
\bigl(W'_k\bigr)\bigr)^2
\bigl(W''_k-{\mathbb E}
\bigl(W''_k\bigr)\bigr)^2 \bigr) + c
\frac{k^5}{n^2}
\nonumber
\end{eqnarray}
for a constant $c <\infty$. Moreover, using the formulas from Lemma
\ref{momentsofVandW},
\begin{eqnarray*}
&& {\mathbb E} \bigl(\bigl(W'_k - {\mathbb E}
\bigl(W'_k\bigr)\bigr)^4 \bigr)
\\
&&\qquad = {\mathbb E}\bigl(\bigl(W'_k\bigr)^4
\bigr) -4 {\mathbb E}\bigl(\bigl(W'_k
\bigr)^3\bigr){\mathbb E}\bigl(W'_k\bigr)+6{
\mathbb E}\bigl(\bigl(W'_k\bigr)^2\bigr)
\bigl({\mathbb E}\bigl(W'_k\bigr)\bigr)^2
\\
&&\quad\qquad{}  -4{\mathbb E}\bigl(W'_k\bigr)
\bigl({\mathbb E}\bigl(W'_k\bigr)\bigr)^3+
\bigl({\mathbb E}\bigl(W'_k\bigr)\bigr)^4
\\
&&\qquad = (1-4+6-4+1) \biggl(\frac{n-k}{n} \biggr)^{i-1}\cdot
\frac{k^2}n +O \biggl( \biggl(\frac{k^2}{n} \biggr)^3+
\biggl(\frac{k^2}{n} \biggr)^2+\frac
{k^2}{n} \biggr).
\end{eqnarray*}
The leading terms cancel, and since $ (\frac{k^2}{n} )^2$ is
dominated by either $ (\frac{k^2}{n} )^3$ or $\frac{k^2}{n}$
depending on $k \geq\sqrt n$ or $k <\sqrt n$, we obtain that
%
%
\begin{equation}
\label{centeredfourthmoment} {\mathbb E} \bigl(\bigl(W'_k - {\mathbb
E}\bigl(W'_k\bigr)\bigr)^4 \bigr) = O
\biggl( \frac
{k^6}{n^3}+\frac{k^2}{n} \biggr)
\end{equation}
and similarly for $W''_k$
%
%
\begin{equation}
\label{centeredfourthmoment1} {\mathbb E} \bigl(\bigl(W''_k
- {\mathbb E}\bigl(W''_k\bigr)
\bigr)^4 \bigr) = O \biggl( \frac
{k^6}{n^3}+\frac{k^2}{n}
\biggr).
\end{equation}
Using the Cauchy--Schwarz inequality and (\ref
{centeredfourthmoment}) and (\ref{centeredfourthmoment1}) in (\ref
{WW}), we get that
\[
\mathbb{V}\bigl(W_k(i)W_k(j)\bigr) \leq c \biggl(
\frac{k^6}{n^3}+\frac
{k^5}{n^2}+\frac{k^2}{n} \biggr)
\]
and therefore
\[
\mathbb{V} (Z_k) \leq c \biggl( \frac{k^5}{n^2}+
\frac{k^2}{n} \biggr)
\]
for some constant $c<\infty$. Plugging this into (\ref{variances}) we
obtain that
%
\begin{eqnarray}
\label{variances2}
\mathbb{V}  ( W_{k-1} - \widetilde W_{k-1} )
&\leq&\biggl(1-\frac{4}{k} \biggr) \mathbb{V} ( W_{k} -
\widetilde W_{k} )\nonumber
\\
&&{}  + c \bigl(\mathbb{V} (W_{k}- \widetilde
W_{k} ) \bigr)^{\nicefrac{1}{2}}\cdot\biggl( \frac{\sqrt k}{n}+
\frac
{1}{k\sqrt n} \biggr)
\\
&&{} + c \biggl(\frac{k}{a_n \sqrt{n}}+\frac{a_n
k}{n^2}+
\frac{1} {k} \biggr).
\nonumber
\end{eqnarray}
Observe that
\begin{eqnarray*}
&& c \bigl(\mathbb{V} (W_{k}- \widetilde W_{k} )
\bigr)^{\nicefrac{1}{2}} \cdot\biggl( \frac{\sqrt k}{n}+\frac{1}{k\sqrt
n}\biggr)
\\
&&\qquad \leq \cases{
\displaystyle \frac{ \mathbb{V} ( W_{k} - \widetilde W_{k} ) }{k},
&\quad if $\displaystyle \bigl(\mathbb{V}(W_{k}- \widetilde W_{k} ) \bigr)^{\nicefrac{1}{2}} \geq c
\biggl( \frac{k^{\nicefrac
{3}{2}}}{n}+\frac{1}{\sqrt n} \biggr) $,
\vspace*{5pt}\cr
\displaystyle 2c^2 \biggl(\frac{ k^2}{n^2}+\frac{1}{kn} \biggr), &\quad else }
\end{eqnarray*}
and therefore since $k \leq a_n$, (\ref{variances2}) becomes
\[
\mathbb{V} ( W_{k-1} - \widetilde W_{k-1} ) \leq\biggl(1-
\frac{3}{k} \biggr) \mathbb{V} ( W_{k} - \widetilde
W_{k} ) + c \biggl(\frac{k}{a_n \sqrt{n}}+\frac{a_n k}{n^2}+
\frac{1} {
k} \biggr).
\nonumber
\]
We now divide both sides by ${k-1\choose 3}$ and iterate up to $a_n$.
Since $\mathbb{V} ( W_{a_n} - \widetilde W_{a_n} ) \leq
c\frac{a_n^2}{n}$, we obtain by Lemma~\ref{momentsofVandW}
\begin{eqnarray*}
\frac{1}{{k-1\choose 3}} \mathbb{V} ( W_{k-1} - \widetilde W_{k-1}
) &\leq& c \cdot\sum_{j=k}^{a_n}
\frac{1}{{j\choose 3}} \biggl( \frac{j}{a_n \sqrt{n} }+ \frac{a_n j} {n^2}+
\frac{1} {j} \biggr)
\\
&&{} + c \cdot\frac{1}{{a_n\choose 3}}\cdot\frac{a_n^2}{n}
\end{eqnarray*}
and therefore
\begin{eqnarray*}
\mathbb{V} ( W_{k-1} - \widetilde W_{k-1} )& \leq& c \cdot
k^3 \biggl( \frac{1}{ka_n \sqrt{n} }+ \frac{a_n} {k n^2}+
\frac
{1}{k^3}+\frac{1}{a_n n} \biggr)
\\
& \leq& c \cdot\biggl( \frac{k^2}{a_n\sqrt{n} } + \frac{a_n k^2} {
n^2} +
\frac{k^3}{a_n n}+1 \biggr).
\nonumber
\end{eqnarray*}
This is the claim.
\end{pf}

\section{Proof of the theorem}\label{sec4}

The proof that $\mu_r$, the expected length of order $r$, is
equal to $\frac{2}{r}$ for every $r \geq1$ can be found in \cite
{Be09}, Theorem~2.11 or in \cite{Du08}, Theorem~2.1. Another quick way
to see this is by using Lemma~\ref{momentsofVandW},
\begin{eqnarray*}
{\mathbb E}\bigl({\mathcal L}^{n,r}\bigr) &=& {\mathbb E} \Biggl( \sum
_{k=2}^n W_k(r)
X_k \Biggr)
\\
& = &\sum_{k=2}^n {
\mathbb E}\bigl(W_k(r)\bigr) {\mathbb E}(X_k)
\\
&=& \sum_{k=2}^n \frac{(n-k) \cdots(n-k-r+2)}{(n-1) \cdots
(n-r)}\cdot
k(k-1) \cdot\frac{1}{{k\choose 2}}
\\
&=& \frac{2}{(n-1) \cdots(n-r)}\cdot\sum_{j=1}^{n-r}
j(j+1) \cdots(j+r-2).
\end{eqnarray*}
The\vspace*{2pt} claim follows now from the fact that $\sum_{j=1}^{n} j(j+1) \cdots
(j+i)=\frac{1}{i+2}n(n+1) \cdots(n+i+1)$.
The asymptotic normality of the total external branch length of
the Kingman coalescent (case $s=1$) was proved in \cite{JK11}. We will
prove the theorem for $s \geq2$.

For $1 \leq r\leq s$ we divide $L^{n,r}$ and the corresponding coupled
quantity into parts. For $1 \leq b_n <a_n\leq n$, let
%
%
\begin{equation}
\label{lengths} \qquad L^{n,r}_{a_n,b_n}:= \sum
_{b_n< k \leq a_n } \frac{2}{k(k-1)} \cdot W_k \quad
\mbox{and}\quad\widetilde{L}^{n,r}_{a_n,b_n}:= \sum
_{b_n<
k \leq a_n } \frac{2}{k(k-1)} \cdot\widetilde{W}_k
\end{equation}
be the length of order $r$ collected between the levels $b_n$ and $a_n$
in the coalescent tree and the corresponding quantity obtained from the
coupling. Note that $L^{n,r}_{n,1}= L^{n,r}$
with $L^{n,r}$ defined in (\ref{L}), and let similarly
%
%
\begin{equation}
\label{Ltilde} \widetilde{L}^{n,r}:= \widetilde{L}^{n,r}_{n,1}.
\end{equation}

Using the coupling we will show that for $\varepsilon>0$
%
\begin{eqnarray}
\label{convprobvector}
&& {\mathbb P} \biggl( \sqrt{\frac{n}{\log
n}}\cdot\bigl\| \bigl(
L^{n,1} - {\mathbb E} \bigl(L^{n,1} \bigr), \ldots,
L^{n,s}- {\mathbb E} \bigl(L^{n,s} \bigr) \bigr)\nonumber
\\
&&\hspace*{56pt}{} - \bigl( \widetilde L^{n,1} - {\mathbb E} \bigl(
\widetilde{L}^{n,1} \bigr), \ldots, \widetilde L^{n,s}- {\mathbb
E} \bigl(\widetilde{L}^{n,s} \bigr) \bigr) \bigr\| \geq{\varepsilon} \biggr)
\to0
\\
\eqntext{\mbox{as } n \to\infty.}
\end{eqnarray}
Once (\ref{convprobvector}) has been proved, the claim of the
theorem follows since the components of the second vector above are by
construction independent and identically distributed and they converge
weakly to the standard normal distribution as $ n \to\infty$, as
follows from the case $s=1$ proved in \cite{JK11}.

The convergence in (\ref{convprobvector}) is a direct consequence of
the following result.
%

\begin{prop}\label{propconvonedim}
For $L^{n,r}$ and $ \widetilde{L}^{n,r}$ defined in (\ref{L}) and
(\ref{Ltilde}), respectively, one has for all $1 \leq r\leq s$ and
$\varepsilon>0$,
\[
{\mathbb P} \biggl( \biggl| \sqrt{\frac{n}{\log n}}\cdot\bigl( L^{n,r} -
{\mathbb E} \bigl(L^{n,r} \bigr) \bigr) - \bigl( \widetilde
L^{n,r} - {\mathbb E} \bigl(\widetilde{L}^{n,r} \bigr) \bigr) \biggr|
\geq{\varepsilon} \biggr) \to0,\qquad\mbox{as } n \to\infty.
\]
\end{prop}

\begin{pf}
We have by the Cauchy--Schwarz inequality that
\begin{eqnarray*}
&& \mathbb{V}   \Biggl(\sum_{b_n< k \leq a_n} \frac{2}{k(k-1)} \cdot
\bigl(W_k -\widetilde{W}_k - \bigl({\mathbb
E}(W_k) - {\mathbb E}(\widetilde{W}_k) \bigr) \bigr)\Biggr)
\\
&&\qquad = \sum_{b_n< k \leq a_n} \sum_{b_n< l \leq a_n}
\frac
{2}{k(k-1)}\frac{2}{l(l-1)} \cdot\mathbb{COV} ( W_{k} -
\widetilde W_{k}, W_{l} - \widetilde W_{l} )
\\
&&\qquad  \leq\sum_{b_n< k \leq a_n} \sum_{b_n< l \leq a_n}
\frac
{2}{k(k-1)}\frac{2}{l(l-1)} \cdot\mathbb{V} ( W_{k} -
\widetilde W_{k} )^{\nicefrac{1}{2}}\mathbb{V} (W_{l} -
\widetilde W_{l} )^{\nicefrac{1}{2}}
\nonumber
\\
&&\qquad = \biggl(\sum_{b_n< k \leq a_n} \frac{2}{k(k-1)} \cdot
\mathbb{V} ( W_{k} - \widetilde W_{k} )^{\nicefrac{1}{2}}
\biggr)^2.
\nonumber
\end{eqnarray*}
Using Lemma~\ref{lemmavar} we obtain
%
\begin{eqnarray}
\label{variancesum}
&& \mathbb{V}
\Biggl(\sum_{b_n< k \leq a_n}
\frac{2}{k(k-1)} \cdot\bigl(W_k -\widetilde{W}_k -
\bigl({\mathbb E}(W_k) - {\mathbb E}(\widetilde{W}_k)
\bigr) \bigr)\Biggr)\nonumber
\\
&&\qquad  \leq c \biggl(\sum_{b_n< k \leq a_n} \frac{1}{k^2}
\cdot\biggl(\frac
{k}{\sqrt{a_n}n^{\nicefrac{1}{4}} } + \frac{\sqrt{a_n}k} {n} + \frac
{k^{\nicefrac{3}{2}} }{\sqrt{a_n n}}+1
\biggr) \biggr)^2
\nonumber\\[-8pt]\\[-8pt]
&&\qquad  \leq c \biggl( \frac{1}{\sqrt{a_n}n^{\nicefrac{1}{4}} } \log\frac
{a_n}{b_n} +
\frac{\sqrt{a_n}} {n} \log\frac{a_n}{b_n} + \frac
{1}{\sqrt{n}}+
\frac{1}{b_n} \biggr)^2
\nonumber
\\
&&\qquad  \leq c \biggl( \biggl( \frac{1}{a_n \sqrt{n}} + \frac{a_n}{n^2} \biggr)
\log^2\frac{a_n}{b_n}+\frac{1}{n}+ \frac{1}{b_n^2}
\biggr).
\nonumber
\end{eqnarray}

In order to show that the claim holds, we consider three regions in the
coalescent tree, namely between level $n$ and level $\frac{n}{(\log
n)^2}$, between level $\frac{n}{(\log n)^2}$ and level $n^{\nicefrac
{1}{2}}$, and finally between level $n^{\nicefrac{1}{2}}$ and level 1,
and write the lengths $L^{n,r}_{n,1}$ and $ \widetilde{L}^{n,r}_{n,1}$
as sums of the lengths gathered in these three regions.

For the first region, let
\[
a_n=n\quad\mbox{and}\quad b_n=\frac{n}{(\log n)^2}.
\]
We obtain from (\ref{variancesum}) and Chebyshev's inequality that
\[
L^{n,r}_{a_n,b_n} - {\mathbb E} \bigl(L^{n,r}_{a_n,b_n}
\bigr) = \widetilde{L}^{n,r}_{a_n,b_n} - {\mathbb E} \bigl(
\widetilde{L}^{n,r}_{a_n,b_n} \bigr) + O_P \biggl(
\frac{\log\log n}{ \sqrt{n}} \biggr)
\]
and therefore
%
%
\begin{equation}
\label{firstregion} \sqrt{\frac{n}{\log n}} \bigl( \bigl( L^{n,r}_{a_n,b_n}
- {\mathbb E} \bigl(L^{n,r}_{a_n,b_n} \bigr) \bigr) - \bigl(
\widetilde{L}^{n,r}_{a_n,b_n} - {\mathbb E} \bigl(
\widetilde{L}^{n,r}_{a_n,b_n} \bigr) \bigr) \bigr) \to0
\end{equation}
in probability as $n \to\infty$.

The second region we consider is the one between the levels $a_n$ and
$b_n$ with
\[
a_n=\frac{n}{(\log n)^2}\quad\mbox{and}\quad b_n=n^{\nicefrac{1}{2}}.
\]
We put together the coupling for the two regions by taking the starting
distribution for the second region to be the distribution of the chain
at the end of the first region. Again from (\ref{variancesum}) we get that
\[
L^{n,r}_{a_n,b_n} - {\mathbb E} \bigl(L^{n,r}_{a_n,b_n}
\bigr) = \widetilde{L}^{n,r}_{a_n,b_n} - {\mathbb E} \bigl(
\widetilde{L}^{n,r}_{a_n,b_n} \bigr) + O_P \biggl(
\frac{1}{\sqrt{n}} \biggr)
\]
and therefore as $n \to\infty$ in probability
%
%
\begin{equation}
\label{secondregion} \sqrt{\frac{n}{\log n}} \bigl( \bigl(L^{n,r}_{a_n,b_n}
- {\mathbb E} \bigl(L^{n,r}_{a_n,b_n} \bigr) \bigr) - \bigl(
\widetilde{L}^{n,r}_{a_n,b_n} - {\mathbb E} \bigl(
\widetilde{L}^{n,r}_{a_n,b_n} \bigr) \bigr) \bigr) \to0.
\end{equation}

For the region in the coalescent between the levels $n^{\nicefrac
{1}{2}}$ and 1, we claim that
%
%
\begin{equation}
\label{E0int} {\mathbb E} \biggl(\sqrt{\frac{n}{\log n}} \cdot
L^{n,r}_{n^{\nicefrac
{1}{2}},1} \biggr) \to0
\end{equation}
and
%
%
\begin{equation}
\label{E0ext} {\mathbb E} \biggl(\sqrt{\frac{n}{\log n}} \cdot\widetilde{L}
^{n,r}_{ n^{\nicefrac{1}{2}},1} \biggr) \to0
\end{equation}
as $n \to\infty$. The second claim follows directly from Proposition~3 in \cite{JK11}, whereas for (\ref{E0int}) we get similarly using
Lemma~\ref{momentsofVandW} that
%
\begin{eqnarray}
\label{expectationlength} {\mathbb E} \bigl(L^{n,r}_{a_n, b_n} \bigr)
&=& {
\mathbb E} \biggl(\sum_{b_n<k \leq a_n} W_k \cdot
\frac{2}{k(k-1)} \biggr)
\\
& \leq& c \sum_{b_n<k \leq a_n} \frac{k^2}{n} \cdot
\frac{1}{k(k-1)} \leq c \cdot\frac{a_n}{n},
\end{eqnarray}
for some constant $c< \infty$. Therefore, setting $a_n=n^{\nicefrac
{1}{2}}$ and $b_n=1$ in (\ref{expectationlength}), we obtain our
claim (\ref{E0int}). Since both $\sqrt{\frac{n}{\log n}} \cdot
L^{n,r}_{ n^{\nicefrac{1}{2}},1}$ and $\sqrt{\frac{n}{\log n}}
\cdot\widetilde{L} ^{n,r}_{ n^{\nicefrac{1}{2}},1}$ are positive
random variables, it follows from (\ref{E0int}) and (\ref{E0ext}),
respectively, that
%
%
\begin{equation}
\label{thirdregion} \sqrt{\frac{n}{\log n}} \cdot L^{n,r}_{n^{\nicefrac
{1}{2}},1}
\to0\quad\mbox{and}\quad\sqrt{\frac{n}{\log n}} \cdot\widetilde{L}
^{n,r}_{n^{\nicefrac{1}{2}},1} \to0
\end{equation}
in probability as $n \to\infty$.

Writing
\[
L^{n,r}=L^{n,r}_{n,1} = L^{n,r}_{n, ({n}/{(\log n)^2})}
+ L^{n,r}_{{n}/{(\log n)^2}, n^{\nicefrac{1}{2}} } +
L^{n,r}_{n^{\nicefrac{1}{2}},1}
\]
and
\[
\widetilde{L}^{n,r}=\widetilde{L}^{n,r}_{n,1} =
\widetilde{L}^{n,r}_{n, ({n}/{(\log n)^2})} + \widetilde
{L}^{n,r}_{{n}/{(\log n)^2}, n^{\nicefrac{1}{2}}}
+ \widetilde{L}^{n,r}_{n^{\nicefrac{1}{2}},1}
\]
and using (\ref{firstregion})--(\ref{E0ext}) and (\ref
{thirdregion}), we get the claim of the proposition, and therefore our
theorem is proved.
\end{pf}

\section*{Acknowledgements} 
We thank two anonymous referees for careful
reading and suggestions that improved the quality of the paper.




\printaddresses

\end{document}